\documentclass[12pt, a4paper]{article}
\makeatletter
\newenvironment{folding}{\endgroup}{\begingroup \def \@currenvir{folding}\edef \@currenvline{\on@line}}
\makeatother
\begin{folding} %Packages
	\usepackage[utf8]{inputenc}
	\usepackage{hyperref }
	\usepackage{latexsym,amssymb, graphicx, bbm, amsmath}
	\usepackage{amsthm, amsfonts, verbatim}
	\usepackage{xcolor}
	\usepackage{pifont}% http://ctan.org/pkg/pifont
	\usepackage{subfig}
	\usepackage{epstopdf}
	
	\AppendGraphicsExtensions{.gif}
	
\end{folding}

\begin{folding} %Tikz
	\usepackage{tikz}
	\usetikzlibrary{shapes.geometric, arrows}
	\tikzstyle{rectwhite} = [rectangle, rounded corners, minimum width=1cm, minimum height=1cm,text centered, draw=black, fill = pink!30]
	\tikzstyle{arrow} = [thick,->,>=stealth]
\end{folding}

\begin{folding} %Charles input
	\usepackage{amsthm, amsfonts, amsmath, amssymb}
	\usepackage{wasysym}
	\usepackage{indentfirst}
	\usepackage{color}
	\usepackage{graphicx}
	\usepackage{microtype}
	\newcommand{\lone}{\lambda_1}
	\newcommand{\ltwo}{\lambda_2}
	\newcommand{\lthree}{\lambda_3}
	\newcommand{\Imag}{\text{Im}}

	\newcommand{\conj}[1]{\overline{#1}}	% conjugado complexo
\end{folding}
\graphicspath{{../figs/}}

\begin{folding} %Theorems
	\newtheoremstyle{dotless}{}{}{}{}{\bfseries}{}{12pt}{}
	\theoremstyle{dotless}
	\newtheorem{thm}{Theorem}[section] % Theorem
	\newtheorem*{thm*}{Theorem} % Theorem (unnumbered)
	\newtheorem*{lem*}{Lemma} % Lemma (unnumbered)
	\newtheorem*{cor*}{Corollary} % Corollary (unnumbered)
	\newtheorem{prop}[thm]{Proposition} % Proposition
	\newtheorem*{prop*}{Proposition} % Proposition (unnumbered)
	\newtheorem{defn}[thm]{Definition} % Definition
	\newtheorem*{defn*}{Definition} % Definition (unnumbered)
	\newtheorem{rem}[thm]{Remark} % Remark
	\newtheorem*{rem*}{Remark} % Remark (unnumbered)
	\newtheorem*{exa*}{Example} % Example (unnumbered)
	\newtheorem*{exe*}{Exercise} % Exercise (unnumbered)
	\newtheorem*{aside*}{Aside} % Aside (unnumbered)
	\newtheorem*{example*}{Example} % Example (unnumbered)
	\newtheorem*{notation*}{Notation} % Notation (unnumbered)
\end{folding}
\begin{folding} %Tikz
	
\end{folding}

\newcommand{\dee}{\mathrm{d}}

\title{
	{Report: Statistics of approximations to zeroes of $\zeta$-function via truncated symmetrized Euler products}\\ \vspace{10mm}
}

\author{Aditya Ghosh}
\date{September 27, 2022}
\begin{document}
	\pagenumbering{gobble}
	\maketitle
	\begin{abstract}
		We look at approximations $ \zeta_X $ of the $\zeta$-function introduced in Gonek's paper \cite{gonek2012finite}. We look at how close the approximate zeroes are to the actual zeroes when (i) X is fixed (Section \ref{sec: stats X5})(ii) X varies like $ t/2\pi $ (Section \ref{sec: stats xt}). We establish a heuristic for estimating these differences, involving values of $ F_X^\star(t) $ and its near-constant slopes near zeta-zero ordinates $ \gamma $. In Section \ref{sec: slopes} we see the slope around the zeroes behaves logarithmically and we calculate a numerical formula for it. In Section \ref{sec: model 1} and \ref{sec: model 2} We scale the differences with the slopes and compare them with models involving 1 or 2 pairs of neighbouring zeta-zeroes. In Section \ref{sec:accuracy of models}, we also look at how often these models capture these scaled differences accurately. In Section \ref{sec: theory}, we look at our methods from a theoretical standpoint. \\ \\
		In Section \ref{sec2: introduction}, we look at how close the approximate zeroes are to the actual zeroes when (i) X is fixed (ii) X varies like $ t/2\pi $. The errors seem to behave like powers of log which should be investigated further from a theoretical standpoint. 
	\end{abstract}
	\pagebreak
	\pagenumbering{arabic}
	\section{Introduction} \label{sec: introduction}
	Our usual convention will be $ s = \sigma + i t $, $ \tau = |t| + 2 $ and $\zeta(s) = \chi(s) \zeta(1-s)$.

	In the critical strip, the values of $ \zeta (s) $ can be well approximated by the \textit{approximate functional equation} \cite{TitchmarshE.C1986Ttot}:
	\begin{equation}\label{eqn: approx func eqn}
		\zeta(s) = \sum_{n \leq X} \frac{1}{n^s} + \chi(s)\sum_{n \leq |t|/2 \pi X} \frac{1}{n^{1-s}} + O(X^{-\sigma}) + O(\tau^{-1/2}X^{1-\sigma})
	\end{equation}
	where $ 0\leq \sigma \leq 1, |s-1| \geq \frac{1}{10}$. To minimize the error terms we take $ X = \sqrt{|t|/2 \pi} $. On the critical line we get an error term of $ O(\tau^{-1/4}) $. 
	\subsection{Approximation by Truncated Symmetrized Euler Product}
	Gonek, Hughes and Keating developed a \textit{hybrid formula} \cite{hybidformula} which aprroximates $\zeta(s)$ in terms of $ P_X(s) $, a weighted version of the truncated Euler Product and $ Z_X(s) $, a expression involving zeta-zeroes close to $ s $:
	\begin{equation}\label{eqn: hybrid formula unmollified}
		\zeta(s) = P_X(s)Z_X(s) + \text{error terms}
	\end{equation}
	where $ P_X(s) := \exp(\sum_{n \leq X} \frac{\Lambda(n)}{n^s \log n}) $. 
	\\
	Using this they conjectured the values of the moments of Zeta Function, by calculating moments of $ P_X $ and $ Z_X $ separately and multiplying them together. These conjectures match the known values of the first few moments. 
	
	Instead of truncating the Dirichlet series and symmetrizing it, we can try to do the same for the Euler product. However, symmetrizing with $ 1 - s $ leads to bad approximations. Instead, Gonek \cite{gonek2012finite} showed symmetrizing with $ \conj{s} $ gives much better results. We define the approximations as: 
	\begin{equation}\label{def: sym trun euler product ZetaX}
		\zeta_X(s) := P_X(s) + \chi(s) P_X(\conj{s}) 
	\end{equation}
	Note that on the critical line, $ \conj{s} = 1-s $. Also, $\zeta_X(s)$ is not an analytic function, instead it is a harmonic function. Assuming RH, it approximates $ \zeta(s) $ well to the right of the critical line. 
	
	It can be easily shown that $ \zeta_X(s) $ satisfies the Riemann Hypothesis: 
	\begin{thm}[\cite{gonek2012finite}, pg.\,2170] 
		Let $ \rho_X = \beta_X + i \gamma_X $ be a zero of $ \zeta_X $ such that $ 0 \leq  \beta_X \leq 1 $ and $ \gamma_X \geq C_0  $. Then $ \beta_X \geq C_0 $, where $ C_0 < 6.3 $ is a constant.
	\end{thm}
	Looking at plots of $ \zeta_X(s) $ and $ \zeta(s) $ on the critical line, we see that even for low values of X, the zeroes of $ \zeta_X $ are quite close to zeta zeroes (figure as provided in \cite{gonek2012finite}):  \\
	\begin{figure}[h]
		\includegraphics[scale = 0.9]{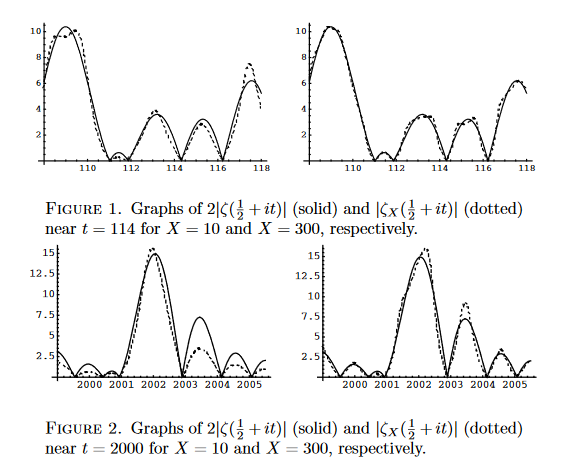}
	\end{figure}

	\subsection{Zeroes of $\zeta_X$}
	From \cite[pg. \,2170]{gonek2012finite} we have that on the critical line, $ |\chi(1/2 + it) | = 1$ for $ |t| \geq C_0 $. Also, $ |P_X(s)| = |P_X(\conj{s})|$. Hence, $ \zeta_X(s) = P_X(s)\left(1 + \chi(s)\frac{P_X(s)}{P_X(s)}\right) $ vanishes for $  |t| \geq C_0 $ if and only if $ \arg \left(\chi(s)P_X(s)/P_X(\conj{s})\right)  \equiv \pi ( \text{mod } 2\pi) $. 
	\begin{prop}[\cite{gonek2012finite}, pg.\, 2171] 
		Define $ F_X(t) = - \arg \chi(1/2 + i t) + 2 \arg P_X(1/2 + it)$. Then for $ |t| \geq C_0 $, $ \zeta_X(1/2 + it) = 0 $ if and only if $ F_X(t) \equiv \pi  ( \text{mod } 2\pi) $
	\end{prop}
	As with Gonek's paper, we work with a slight modification of $ F_X $ to simplify our expressions. We define: 
	\begin{defn}
		\begin{eqnarray*}
			P_X^\star(s) &=& P_X(s) \exp (-F_2((s-1)\log X))\\
			\zeta_X^\star(s) &=& 	P_X^\star(s) + \chi(s)	P_X^\star(\conj{s})\\
			F_X^\star(t) &=& - \arg \chi(1/2 + i t) + 2 \arg P_X^\star(1/2 + it)
		\end{eqnarray*}
	\end{defn}
	
	We then have $ F_X^\star(t) =  F_X(t) + O(\frac{X}{\tau^2 \log X})$, hence for values of $ X \leq t^2 $, these are close. For $  F_X^\star(t) $ we have the simpler expression: 
	\begin{prop}[\cite{gonek2012finite}, pg.\,2180] The zeroes of $ \zeta_X^\star(1/2 + it) $ are solutions of $ F_X^\star(t) \equiv \pi  ( \text{mod } 2\pi) $. Assuming the Riemann Hypothesis, then:
		\begin{equation} \label{eqn: Fxstar explicit formula} 
			\frac{1}{2 \pi}  F_X^\star(t) = N(t) - 1 - \frac{1}{\pi} \Imag  \sum_{\gamma} F_2(i(t-\gamma) \log X) + O(\frac{X^{-3/2}}{\tau^2 \log X}) 
		\end{equation}
		where the sum is over the ordinates $\gamma$ of zeta-zeroes. $ F_2(z) = 2 E_2(2z) - E_2(z)  $, where $ E_2(z) = \int_{z}^{\infty} \frac{e^{-w}}{w^2} \dee w $ is the second exponential integral. $ N(T) $ is the number of zeta zeroes with ordinate in $ [0,T] $. 
	\end{prop}
	We compare the graphs of $ N(t) $ and $ 1 +	\frac{1}{2 \pi}  F_X^\star(t) $:
	\setcounter{figure}{2} 
	\begin{figure}[h]
		\centering
		\includegraphics[scale=0.8]{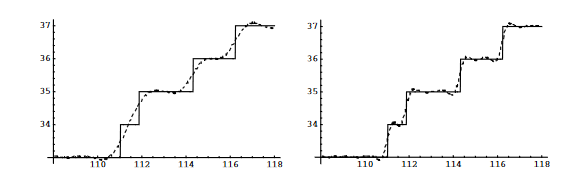}
		\caption{Graphs of $ N(t) $ and $ \dfrac{1}{2 \pi}F_X(t) + 1 $ near $ t = 114 $ for $ X = 10 $ and $ X= 300 $, as provided in \cite{gonek2012finite}}
		\label{fig: staircse}
	\end{figure}
	\subsection{The heuristic idea of the article}
	As we see in the figure \ref{fig: staircse}, the actual zeroes and approximate zeroes are quite close to each other. Also, around a zeta zero the slope stays relatively constant. For the zeta-zero $ \gamma_0 $, let the approximate zero be $ \tilde{\gamma_0} $, and their difference is $ \delta = \gamma_0 - \tilde{\gamma_0} $. 
	
	Let $  H_X^\star:= \frac{1}{2\pi}  F_X^\star + 1 $. At $ \tilde{\gamma_0} $, $ H_X^\star $ takes the value $ N(T) -  \frac{1}{2} $. Then,  by simple geometry we have 
	\begin{equation}
		(\text{Slope around } \gamma_0) \times \delta =  H_X^\star(\gamma_0) - (N(\gamma_0) - \frac{1}{2})
	\end{equation}
	Looking at \ref{eqn: Fxstar explicit formula}, at $ t = \gamma_0 $, 
	\begin{equation}
		H_X^\star(\gamma_0) = N(\gamma_0) - \underbrace{\frac{1}{2}}_{\text{from } \gamma = \gamma_0} - \frac{1}{\pi} \Imag  \sum_{\gamma \neq \gamma_0} F_2(i(\gamma_0-\gamma) \log X)
	\end{equation}
	Hence, we have: 
	\begin{equation} \label{eqn: guiding equation}
		(\text{Slope around } \gamma_0) \times \delta = - \frac{1}{\pi} \Imag  \sum_{\gamma \neq \gamma_0} F_2(i(\gamma_0-\gamma) \log X)
	\end{equation}
	$ F_2 $ decreases for bigger values, according to $ F_2(z) \ll \frac{\exp(\max(- \Re \, z, - \Re\, 2z))}{|z|^2} $. Hence we look at the contribution of zeroes close to $ \gamma_0 $.
	\section{Statistics for zeroes of $\zeta_X$ and $ \zeta $ for X = 5} \label{sec: stats X5}
	\subsection{Actual Statistics} 
	We observe some of the statistics of zeroes when X is fixed, say X=5:
	\newpage
	\begin{figure}[!h]
		\centering
		\subfloat[List Plot]{{\includegraphics[scale = 0.7]{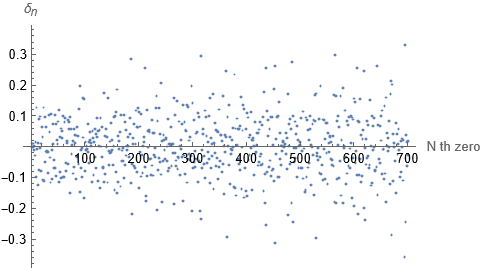} }}
		\qquad
		\subfloat[Histogram]{{\includegraphics[scale = 0.7]{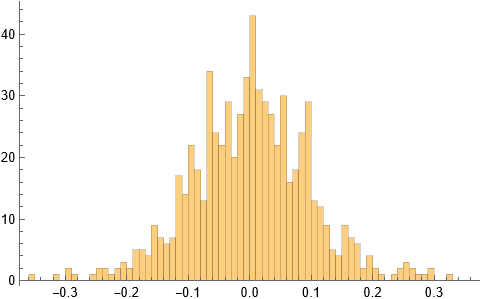} }}
		\caption{Difference between actual and approximate first 700 zeroes for X=5}
		\label{fig: X=5 actual stats}
	\end{figure}
	We can see that most of the differences are very small but in general the differences tend to increase as we look at higher zeroes. 
	
	Looking at the two neighbouring zeroes, we define:
	\begin{equation}\label{eqn: X5 first approx}
		h1(n) = - \frac{1}{\pi} \Imag F_2(i  (\gamma_n - \gamma_{n-1}) \log 5) +  \frac{1}{\pi} \Imag F_2(i  (\gamma_{n+1} - \gamma_{n}) \log 5) 
	\end{equation}
	where $ \gamma_{n} $ is the $ n^\text{th} $ zeta-zero ordinate.
	\subsection{Model Statistics}
	Similarly we can look at more neighbouring zeroes. However this approach is possibly fruitless as we consider higher ordinates, justified in Section \ref{sec: bad idea}. 
	
	\section{Eigenvalue spacings of $ 3 \times 3  $ GUE}
	Looking at \ref{eqn: Fxstar explicit formula}, the $ \log X  $ factor in the parenthesis is like a scaling factor. To apply Random Matrix Theory, it would be helpful to scale the zeroes appropriately. To that end, to approximate the expression around $ t = \gamma_0 $, we choose $ X = {\gamma_0}/{2 \pi} $. Thus, for zeta-zero ordinates $ \gamma $ close to $ \gamma_0 $, the scaling by $ \log X = \log \frac{\tau}{2\pi}  $ ensures they behave like scaled eigenvalues of a GUE.
	
	We quickly derive the joint probability distribution of eigenangle spacings of a $ 3 \times 3 $ GUE. Given a random GUE matrix, label the eigenvalues  $\lambda_1\geq\lambda_2\geq\lambda_3$. Then, $ {P}(\lambda_1,\lambda_2,\lambda_3) = C e^{-\frac{1}{2}(\lambda_1^2 + \lambda_2^2 + \lambda_3^2)}(\lambda_1 - \lambda_2)^2(\ltwo -\lthree)^2(\lone-\lthree)^2$ as usual. 
	\begin{align*}
		&{P}(\lambda_1 - \lambda_2 = \alpha, \ltwo - \lthree = \beta ) \\
		&=  \int_{\lthree = - \infty}^{\infty} \int_{\ltwo = \lthree}^\infty \int_{\lone = \ltwo}^{\infty} {P}(\lambda_1,\lambda_2,\lambda_3) \delta(\alpha - (\lone - \ltwo)) \delta(\beta - (\ltwo - \lthree)) \, \dee \lone \dee \ltwo \dee \lthree \\
		&= \int_{\lthree = -\infty}^{\infty} P(\alpha + \beta + \lthree, \beta + \lthree, \lthree) \dee \lthree\\
		&= C \int_{\lthree = -\infty}^{\infty} \alpha^2 \beta^2 (\alpha + \beta)\exp(- ((\alpha + \beta + \lthree)^2 + (\beta + \lthree)^2 + \lthree^2)) \dee \lthree\\
		&= C \alpha^2 \beta^2 (\alpha + \beta)^2 e^{-\frac{1}{2}(\alpha + \beta)^2} e^{- \beta^2} \int_{\lthree = -\infty}^{\infty} \exp(- 3 \lthree^2 -2(\alpha + 2 \beta)\lambda_3) \dee \lthree \\
		&=  C \alpha^2 \beta^2 (\alpha + \beta)^2 e^{-\frac{1}{2}(\alpha + \beta)^2} e^{-\frac{1}{2} \beta^2} \int_{\lthree = -\infty}^{\infty} \exp(-\frac{3}{2}(\lthree + \frac{\alpha + 2 \beta}{3})^2) \exp(\frac{(\alpha + 2 \beta)^2}{6}) \dee \lthree \\
		&= \tilde{C} \alpha \beta (\alpha + \beta) e^{-\frac{1}{2}(\alpha + \beta)^2} e^{-\frac{1}{2} \beta^2} e^{\frac{(\alpha + 2 \beta)^2}{6}}\\
		&= \tilde{C} \alpha^2 \beta^2 (\alpha + \beta)^2 e^{-\frac{1}{3}(\alpha^2 + \alpha\beta + \beta^2)}
	\end{align*}

	Scaling them properly by $ \phi_k = \lambda_k  $, we have: 
	\begin{prop} \label{prop: 33 GUE Spacings}
		The joint distribution of scaled eigenvalue spacings of a $ 3 \times 3 $ GUE is: 
		$ {P}(\phi_1 - \phi_2 = \alpha, \phi_2 - \phi_3 = \beta ) =  \tilde{C} \alpha^2 \beta^2 (\alpha + \beta)^2 e^{-\frac{1}{3}(\alpha^2 + \alpha\beta + \beta^2)}$
		Numerically we find $ C^\prime $ to be around 0.046
	\end{prop}
	\section{Zeroes of $ \zeta_X $ with $ X = \lfloor \frac{t}{2 \pi} \rfloor $}
	\subsection{Actual statistics} \label{sec: stats xt}
	As we mentioned before, to approximate the expression around $ t = \gamma_0 $, we choose $ X = \frac{\gamma_0}{2 \pi} $.
	\begin{figure}[!h]
		\centering
		\subfloat[List Plot]{{\includegraphics[scale = 0.75]{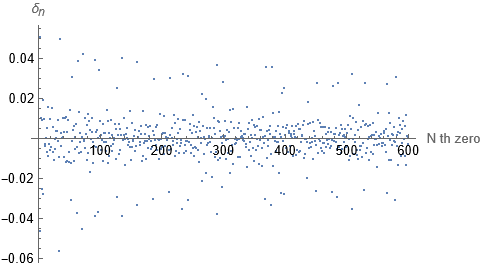} }}
		\qquad
		\subfloat[Histogram]{{\includegraphics[scale = 0.7]{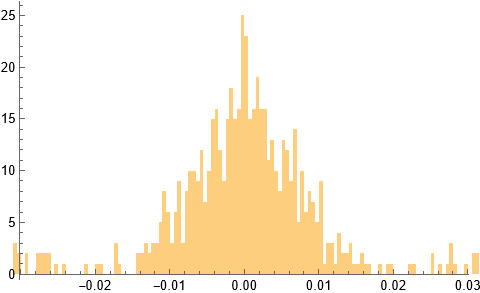} }}
		
		\caption{Difference between actual and approximate first 700 zeroes for  $ X = \frac{\gamma_0}{2 \pi} $}
		\label{fig: Xt Actual stats}
	\end{figure}\\
	We can see that most of the differences are very small but unlike  the differences tend to increase as we look at higher zeroes. 
	\subsection{Scaling with the near-constant slopes around $ \gamma_0 $} \label{sec: slopes}
	Let us first look at some numerical evidence.  
	\begin{figure}[h!]
		\includegraphics[scale = 0.5]{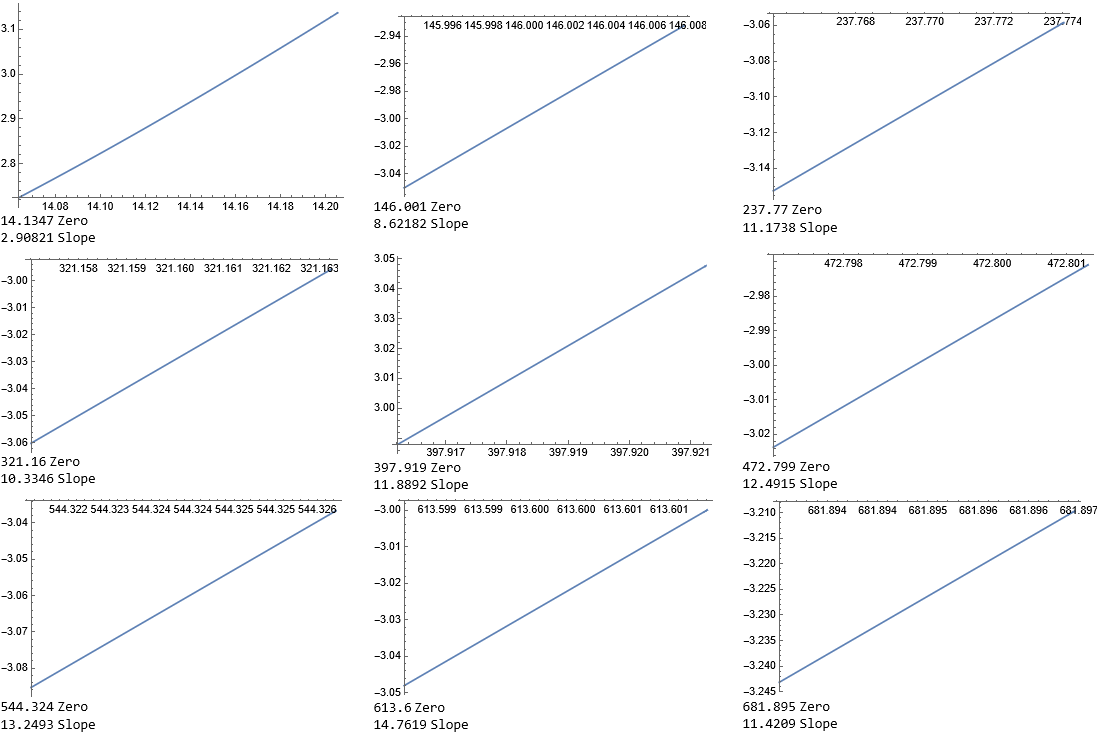}
		\centering
		\caption{Graph of $ F_X^\star $ in interval $ [\gamma - \frac{1}{\gamma}, \gamma + \frac{1}{\gamma}] $ for various zeta-zero ordinates $ \gamma $, along with corresponding slope in that region.}
		\label{fig: Slopes collage}
	\end{figure}
	\\
	As one can see, the slopes are near constant in the neighbourhood of the zeta-zero ordinates.
	\pagebreak
	
	We calculate the slopes and plot them: 
	
	\begin{figure}[h!]
		\includegraphics[scale = 1]{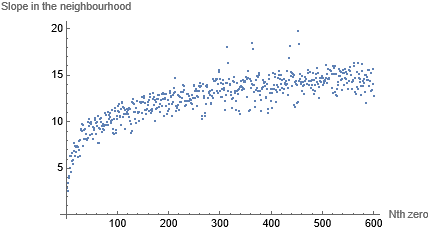}
		\centering
		\caption{Slope of $ F_X^\star $ in interval $ [\gamma - \frac{1}{\gamma}, \gamma + \frac{1}{\gamma}] $ for upto 600 zeta-zero ordinates $ \gamma $}
		\label{fig: SlopeNumerics}
	\end{figure}
	
	The slopes closely resemble a logarithmic graph. We can deduce a formula for it later. For now, we work with the numerical values.
	So, according to our heuristic , scaling the differences by our slopes should give us the right hand side of equation \ref{eqn: guiding equation}. We do so here \footnote{I had to weed out some anomalies (about 10 whose values were way off). There were probably because of errors in calculating our near-constant slope due to shortcomings in our code. I will { \color{red} later try to rescale the differences again with the theoretical value} for the slope.}: 
	
	\begin{figure}[!h]
		\includegraphics[scale = 0.8]{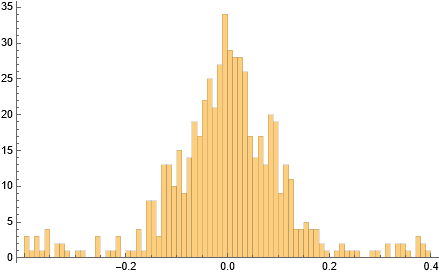}
		\centering
		\caption{Scaled differences between actual approximate first 600 zeroes}
		\label{fig: Xt Scaled Zero diff numerics}
	\end{figure}
	From Gonek's paper \cite[pg./,2175]{gonek2012finite}, we have, $$ F_x^\prime(t) = \log(\frac{t}{2\pi}) - 2 \sum_{n \leq X^2 }\frac{\Lambda_X(n) \cos(t \log n)}{\sqrt{n}} + O(1/\tau) $$
	Changing the x-axis of our graph from the `$ N^{\text{th}} $ zero' to `Ordinate t', we use Mathematica to get the best fit log curve as $ -0.6 + 3.1 \log (\frac{t}{2 \pi}) $. The R-Squared value is 99.327 \%. 
	I've plotted the data with the fitted curve below: 
	\begin{figure}[h!]
		\centering
		\includegraphics[scale = 0.7]{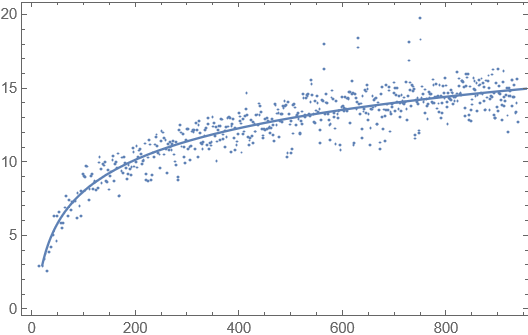}
		\caption{Slope (points) and Log Model (solid line) vs Ordinate $ t $}
		\label{fig: fitted slope}
	\end{figure}
	\subsection{Model statistics : with 1 pair of adjacent zeroes} \label{sec: model 1}
	Similar to \ref{eqn: X5 first approx}, looking at the two neighbouring zeroes, we define:
	\begin{equation}\label{eqn: Xt first approx}
		h1(n) = - \frac{1}{\pi} \Imag F_2(i  (\gamma_n - \gamma_{n-1}) \log \dfrac{\gamma_n}{2 \pi}) +  \frac{1}{\pi} \Imag F_2(i  (\gamma_{n+1} - \gamma_{n}) \log \dfrac{\gamma_n}{2 \pi}) 
	\end{equation}
	where $ \gamma_{n} $ is the $ n^\text{th} $ zeta-zero ordinate.
	
	Plotting it for the first 2000 zeroes, we get:
	\newpage
	\begin{figure}[!h]
		\centering
		\subfloat[List Plot]{{\includegraphics[scale = 0.6]{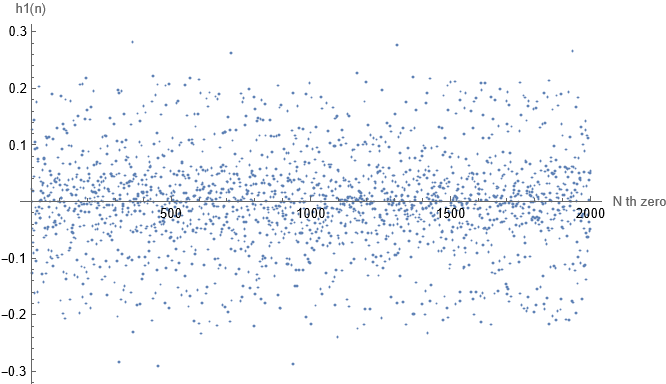} }}
		\qquad
		\subfloat[Histogram]{{\includegraphics[scale = 0.8]{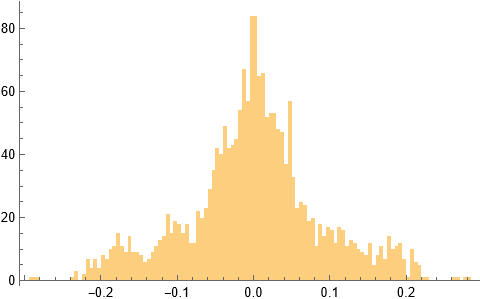} }}
		\caption{Values of model h1(n) for the first 2000 zeroes}
		\label{fig: Xt first model}
	\end{figure}
	
	If we squint a bit, this resembles the shape of the distribution we get in our actual scaled differences in Figure \ref{fig: Xt Scaled Zero diff numerics}. More notably, the majority of the values are recorded between -0.2 and 0.2 in both the figures.
	
	The main difference between the model and the actual histogram is that the model lacks a lot of values which lie in the region 0.05 to 0.15 (and symmetrically on the other side too). It can be heuristically explained by the fact that, in the infinite sum, we can pair up zeroes on the left of $ \gamma_0 $ to the right of it. Most of the time they cancel out as they have similar absolute values. The values in this region 0.05 to 0.15 arises when these pairings don't cancel each other out, which is more likely once we consider more and more pairings instead of just 1 pair. 
	\subsection{Model statistics : with 2 pairs of adjacent zeroes} \label{sec: model 2}
	As a follow-up to our last observation we look at the statistics for two pairs of adjacent zeroes. Defining an analogous model as \ref{eqn: Xt first approx}, we consider two pairs of adjacent zeroes. We place the results one after the other for comparison \footnote{I've cleaned up some of the more extreme values ($ \geq 0.3 $) in \ref{fig: Xt Scaled Zero diff numerics} so that the histograms match up in their width}:
	\begin{figure}[!htb]
		\begin{minipage}{0.32\textwidth}
			\includegraphics[width=\linewidth]{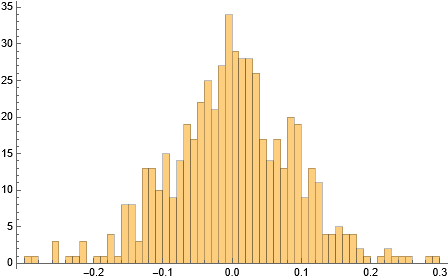}
			\caption{Actual scaled differences (first 600 zeroes) }\label{fig:awesome_image1}
		\end{minipage}\hfill 
		\begin{minipage}{0.32\textwidth}
			\includegraphics[width=\linewidth]{Xtfirstmodel}
			\caption{Model using 1 pair of adjacent zeroes(first 2000 zeroes)}\label{fig:awesome_image2}
		\end{minipage}\hfill
		\begin{minipage}{0.32\textwidth}%
			\includegraphics[width=\linewidth]{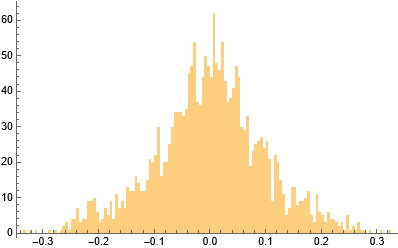}
			\caption{Model using 2 pairs of adjacent zeroes (first 2000 zeroes)}\label{fig:awesome_image3}
		\end{minipage}
	\end{figure}
	
	As we had remarked before, due to more scope of uneven cancellation between pairs of zeroes, there is more  ``meat" around the central spike in Figure \ref{fig:awesome_image3}.
	\subsection{Accuracy of each model}\label{sec:accuracy of models}
	What we want to study in this section is, what is the probability that the first/second model gives answers close to the actual scaled difference. The motivation for this lies in remark \ref{remark1} in the next section.
	
	At a given zeta-zero ordinate $ \gamma_n $,we have the actual scaled difference $ D_n $ and the value given by the model $ M_n $. The relative error in the model be $ E_n = \frac{|D_n - M_n|}{D_n} $. Given a sample of values of n from $ N $ to $ 600 $, we want to calculate for what percentage of the sample $ E_n < \frac{1}{5} $.  We do this for the first and second model:  \newpage
	
	\begin{figure}[!htb]
		\centering
		\begin{minipage}{0.4\textwidth}
			\includegraphics[width=\linewidth]{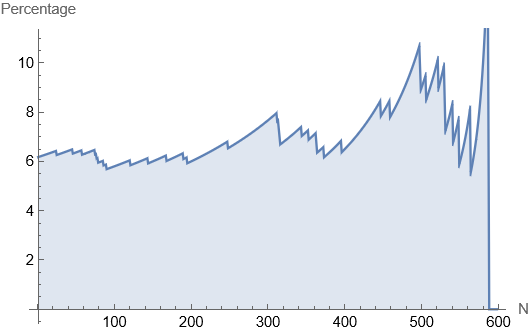}
			\caption{Model 1, Percentage of samples from N th zero to 600th that have \
				error $ <  1/5 $ th the actual value" }\label{fig: percentage model 1}
		\end{minipage} \hspace{2 em}
		\begin{minipage}{0.4\textwidth}
			\includegraphics[width=\linewidth]{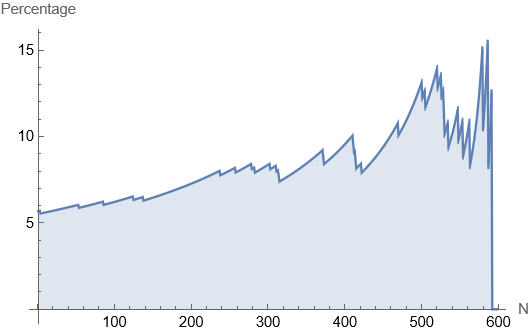}
			\caption{Model 2, Percentage of samples from N th zero to 600th that have \
				error $ <  1/5 $ th the actual value" }\label{fig: percentage model 2}
		\end{minipage}
	\end{figure}

	Ideally we are looking for large values of ordinates. We see that for the first model, we have small error for $ 6\% $ of the time, whereas for the second model it is around $ 10\% $.   
	\section{Theoretical Justifications} \label{sec: theory}
	We begin with recalling some properties from Gonek's paper \cite{gonek2012finite}:
	\begin{enumerate}
		\item (pg.\,2164) $ F_2(z) \ll \frac{\exp(\max(- \Re \, z, - \Re\, 2z))}{|z|^2} $
		\item (Lem 3.3) Assuming RH, $ \sum_{|\gamma - t| > \Delta} \frac{1}{(t-\gamma)^2}\ll \frac{1}{\Delta}\left(\log \tau + \frac{\Phi(\tau)}{\Delta}\right)$ where $ \Phi(t) $ is admissible 
		\item $ \Phi(t) = \dfrac{1}{6} \log \tau $ is admissible
	\end{enumerate}
	We split up the sum in \ref{eqn: Fxstar explicit formula}, 
	\begin{equation}
		\begin{aligned}
			& \Imag \sum_{\gamma \neq \gamma_0} F_2(i (\gamma_0 - \gamma) \log X) \\
			& \qquad=\Imag \sum_{ 0< |\gamma - \gamma_0| \leq \Delta} F_2(i (\gamma_0 - \gamma) \log X) + \Imag \sum_{|\gamma - \gamma_0| > \Delta} F_2(i (\gamma_0 - \gamma) \log X) 
		\end{aligned}
		\label{eqn: Approx F_X with tails}
	\end{equation}
	\subsection{Why keeping X fixed is a bad idea} \label{sec: bad idea}
	\noindent Let X be fixed. If we try to bound the second term, we have: 
	\begin{align*}
		\Imag \sum_{|\gamma - \gamma_0| > \Delta} F_2(i (\gamma_0 - \gamma) \log X) &\ll 	\Imag \sum_{|\gamma - \gamma_0| > \Delta} \frac{1}{|\gamma_0 - \gamma|^2 \log^2 X}\\
		&\ll \frac{1}{\Delta \log^2 X}\left(\log \tau + \frac{\Phi(\tau)}{\Delta}\right)
	\end{align*}
	Since, $ X $ is fixed, to make this term negligible, we must take $ \Delta$ large such that $ \log \tau = o(\Delta) $. However, the first term would then go over zeroes with $ 0 < |\gamma_0 - \gamma| \leq \Delta $ which would include a lot of neighbouring zeroes for high ordinates. Hence, approximating by the neighbouring two zeroes is pointless. 
	\subsection{Why varying X like $ t/{2\pi} $ is a good idea}
	Looking at $ X = \lfloor\gamma_0/{2 \pi} \rfloor \sim \tau/{2 \pi} $, the second term in Equation \ref{eqn: Approx F_X with tails} can be bounded as: 
	\begin{equation} \label{eqn: tail Xt}
		\begin{aligned}
			\Imag \sum_{|\gamma - \gamma_0| > \Delta} F_2(i (\gamma_0 - \gamma) \log X) &\ll 	\Imag \sum_{|\gamma - \gamma_0| > \Delta} \frac{1}{|\gamma_0 - \gamma|^2 \log^2 X}\\
			&\ll \frac{1}{\Delta \log^2 X}\left(\log \tau + \frac{\Phi(\tau)}{\Delta}\right) \vspace{1.5 em}
			&{\small \text{(plug in $ \Phi(t) = 1/6 \log \tau $)}}\\
			& \sim \frac{1}{\Delta} \left(\frac{1}{\log \tau} + \frac{1}{6 \Delta \log \tau }\right)
		\end{aligned}
	\end{equation}
	If we take $ \Delta = C \frac{2 \pi}{\log \tau} $. We take constant C such that for large ordinates this term has a smaller order than the first term .
	
	Now, looking at the first term we have: 
	\begin{align*}
		\Imag \sum_{ 0< |\gamma - \gamma_0| \leq \Delta} F_2(i (\gamma_0 - \gamma) \log X) &= \Imag \sum_{ 0< |\gamma - \gamma_0| \leq \Delta} F_2(2\pi i (\gamma_0 - \gamma) \frac{\log (\gamma_0/2\pi)}{2 \pi })\\
		&= \Imag \sum_{ 0< |\gamma - \gamma_0| \frac{\log \gamma_0 }{2 \pi}\leq C} F_2(2 \pi i (\gamma_0 - \gamma) \frac{\log (\gamma_0/2\pi)}{2 \pi })
	\end{align*}
	\begin{rem} \label{remark1}
		Around the point $ \gamma_0 $, the scaled ordinates of zeta-zeroes behave like eigenvalues of $ N \times N $ GUE, where $ N = \lfloor\log \gamma_0\rfloor $. Given a suitable C, we can calculate the probability with which there will be only 1 pair of adjacent scaled eigenvalues at a distance C. This is the probability that our model with 1 pair of adjacent zeroes approximates our actual statistics well. Similarly we can do this for k-pairs of adjacent zeroes.
	\end{rem}
	\begin{rem}
		If we have a better estimate of this tail in \ref{eqn: tail Xt}, we can make do with much smaller values of $\Delta$ which will remove a lot of these complications. 
	\end{rem}

	\section{Experimental Values for higher ordinates} \label{sec2: introduction}
	In the previous sections we looked at two cases:
	\begin{enumerate}
		\item X fixed at $ X=5 $ for up to first 700 zeroes 
		\item X varied as $ X = \frac{t}{2\pi} $ for up to first 600 zeroes
	\end{enumerate}
	The limitations were primarily due to the function $ \chi(s) = \pi^{s-\frac{1}{2}} \frac{\Gamma(\frac{1-s}{2})}{\Gamma(\frac{s}{2})} $ used to define the approximate zeta function $ \zeta_X $. Mathematica could not handle values as small as $ e^{-700} $ which would arise from these .We replace $ \chi(s) $ with an estimate arising from Stirling's approximation: 
	\begin{equation}\label{chiapp}
		\chi(s) = \left(\frac{\tau}{2 \pi }\right)^{\frac{1}{2}-\sigma- i t} e^{it + \frac{1}{4}i \pi } \left(1 + O(1/\tau)\right)
	\end{equation}
	This allows us to expand the scope of our investigations drastically. I have computed the following things:
	\begin{enumerate}
		\item X fixed at $ X = 3$ to $27 $ and$  X = 30,35,\dots,65 $ and $  X = 100 $ for up to first 20,000 zeroes. 
		\item  X varied as $ X = \frac{t}{2\pi} $ for up to first 5000 zeroes and zeroes 10,000 to 11,000
		\item X fixed at $ X = 20 $ for up to first 100,000 zeroes
	\end{enumerate}
	
	We calculate the variances of the errors and scale them by the density at that point, that is,   $V_n = \left({\text{(Error at n-th zero)} \times  \frac{1}{2\pi} \log(\frac{\gamma_n}{2\pi})}\right)^2 $. To make things a bit smoother we look at the moving averages of these $ V_n $, averaged over the neighbouring 100 zeroes.  
	
	{\color{red} Something to note here is that unlike before we don't scale the differences by the slope of $ F_X(t) $ function which is around $ 3 \log(\frac{t}{2 \pi}) $. Instead we scale by $ \frac{1}{2\pi}\log(\frac{t}{2\pi}) $ which is the inverse of the density of zeta-zeroes at that point. It only changes the values by a constant factor.}
	\subsection{Moving Scaled Variances for X fixed}
	We calculate the variances for $ X = 3$ to $27 $ and $  X = 30,35,\dots,65 $ and $  X = 100 $. Here's what they look like for some of them: 
	\begin{figure}[h]
		\centering
		\includegraphics[scale=0.5]{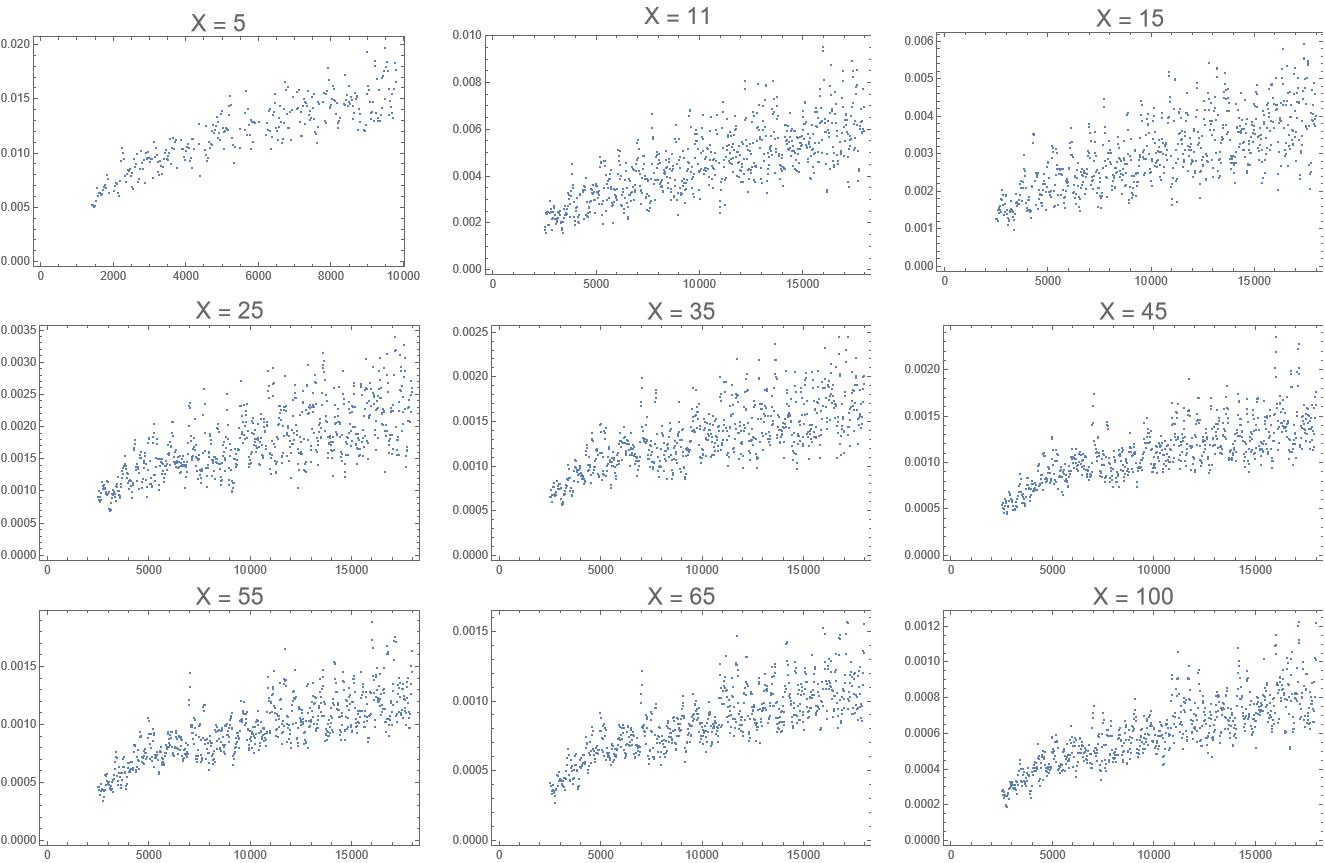}
		\caption{Plot of Moving variance vs Ordinate of zeta-zero }
		\label{fig: plotmovingvartable}
	\end{figure}

	As we can see the plots resemble some sort of logarithmic growth as the ordinate increases. Also notice that as X increases the variances also decrease. This is expected. 
	\subsection{Scaled variances for X varied like $X = \frac{t}{2 \pi}$}
	We calculate the absolute errors as X varies like $X = \frac{t}{2 \pi}$ and then scale them by inverse of mean density as before. We do this for up to first 5000 zeroes and zeroes 10,000 to 11,000. We take the moving variance (200 points) of these values to make it more smooth. The results are: 

	\begin{figure}[h]
		\centering
		\includegraphics[scale=0.8]{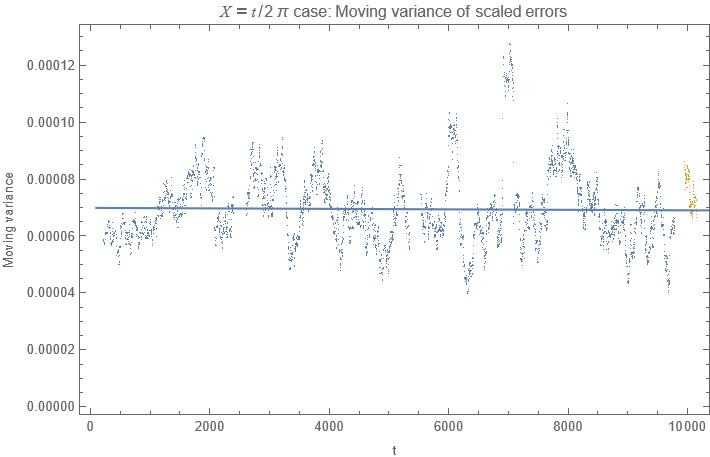}
		\caption{Moving variance for $ X = \frac{t}{2\pi} $ and best fit line $ V(t) = 6.9865 \times 10^{-5} - 8.54 \times 10^{-11} t $}
		\label{fig: plotxvar}
	\end{figure}
	In the previous report we had hoped that these errors would remain bounded uniformly even for high ordinates and that looking at a few adjacent zeroes would give a good idea of the distribution of errors. Figure \ref{fig: plotxvar} confirms this.
	
	If we look at the histogram of scaled errors at various points of the sample, we see the following trend: \newpage
	
	\begin{figure}[h]
		\centering
		\includegraphics[scale=0.7]{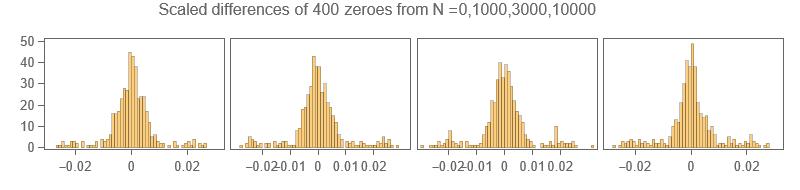}
		\caption{Histogram of scaled errors (not variances!)}
		\label{fig: histerror}
	\end{figure}

	The distribution remains similar for higher ordinates.
	
	We recall what we observed in the first report where we took more and more adjacent zeroes in our models\footnote{Another reminder that here we are scaling by the slope and not by the density, but it only amounts to a factor of $ \frac{3}{2 \pi} $}: \\ 
	
	\begin{figure}[h]
		\begin{minipage}{0.32\textwidth}
			\includegraphics[width=\linewidth]{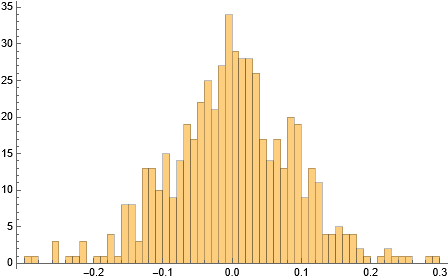}
			\caption{Actual scaled differences (first 600 zeroes) }\label{fig:awesome_image1}
		\end{minipage}\hfill 
		\begin{minipage}{0.32\textwidth}
			\includegraphics[width=\linewidth]{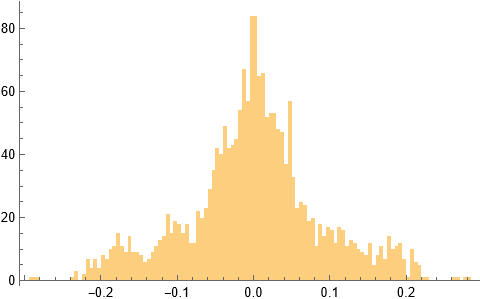}
			\caption{Model using 1 pair of adjacent zeroes(first 2000 zeroes)}\label{fig:awesome_image2}
		\end{minipage}\hfill
		\begin{minipage}{0.32\textwidth}%
			\includegraphics[width=\linewidth]{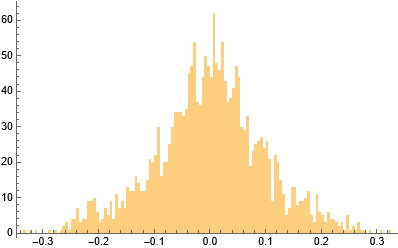}
			\caption{Model using 2 pairs of adjacent zeroes (first 2000 zeroes)}\label{fig:awesome_image3}
		\end{minipage}
	\end{figure}  

	Hence our initial guess that using a few pairs of adjacent zeroes can model the error distribution consistently is true. Otherwise the distribution would've become ``fatter" for higher ordinates to incorporate more uneven cancellations from more terms.
	
	\subsection{Moving variances as  $ V(X,t) = A(X)B(t) $ }
	We try to find a model for the moving variances we had in Figure \ref{fig: plotmovingvartable}. Let's say the moving variances behave like $ A(X)B(t) $. The method I have used is:
	\begin{enumerate}
		\item Guess a model for $ B(t) $ such as $  \log^k(t) $, $  t^k $, $ W(t)^k $ (Lambert W function) or $ t^{k_1} \log^{k_2}(t) $. \\
		For example, take $ B(t) =  \log^k(t)$ here:\\
		\begin{figure}[!h]
			\begin{minipage}{\textwidth}
				\centering
				\includegraphics[scale=0.33]{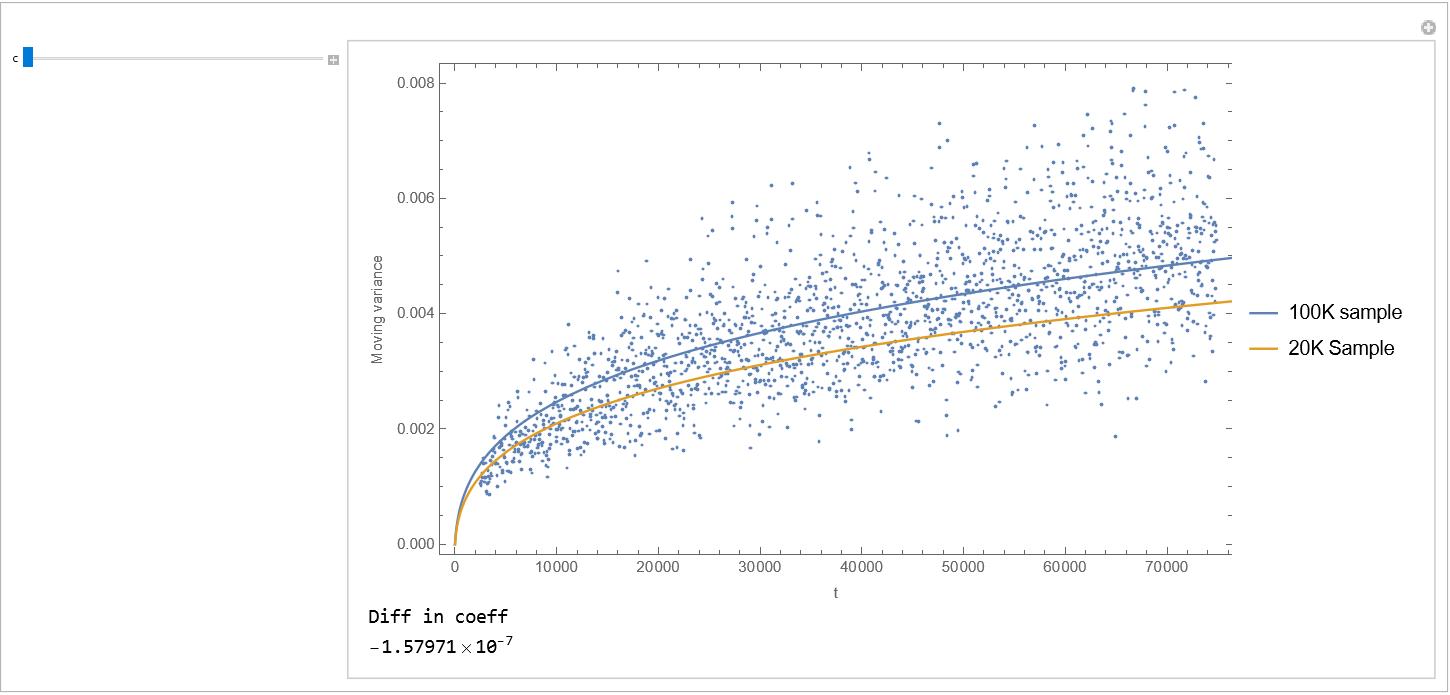}
				\caption{Best fit curves for $ B(t) = \log^k(t) $ with $ k = 3.5 $}
				\label{fig: Bestfitlogpow}
			\end{minipage}
		\end{figure}\\
		\item With the ability to manipulate k, take $ X = 20 $. Remember, we have the moving variance for this upto $ 10^5 $ zeroes. First restrict to the sample of $ 5000 $ zeroes and find the best fit coefficient $ a_1 $ to our chosen model of $ B(t) $. Now take the full sample of $ 10^5 $ zeroes and find the best fit coefficient $ a_2 $. 
		\item Vary k to find the point where the difference $ |a_1 - a_2| $ is the lowest. This is the value of $ k $ we will consider henceforth. 
		
		In the example of $ \log^k(t) $ considered in Figure \ref{fig: Bestfitlogpow}, $ k = 4.3 $ seems to have the lowest difference in coefficients.
		\item Now, vary X and at each step find the best fit coefficient for our model $ B(t) $ (with the value of k found in step 3). These will be the values for A(X). 
		\item Plot the values of $ A(X) $ for the different X and find the best fit model in a similar way to $ B(t) $. 
		
		In our example of $ B(t) = \log^{4.3}(t)$, we find the coefficients A(X) have the following values: \newpage
		
		\begin{figure}[!h]
			
			\centering
			\includegraphics[scale=0.7]{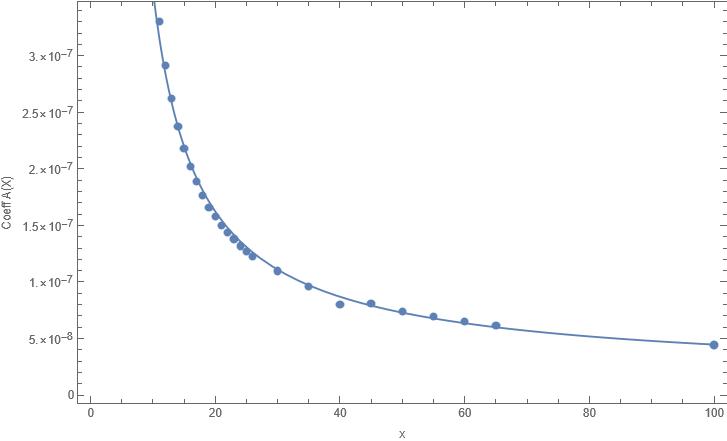}
			\caption{Coefficients $ A(x) $ for $ B(t) = \log^{4.3}(t) $}
			\label{fig: logpowAX}
			
		\end{figure}
	
		Here, $ A(X) = \frac{5.02365 \times 10^{-6}}{\log^{3.13788}(X)} $ is the best fit curve derived from sample X = 11 to 26, but also seem to model the values from X = 30 to 100 very well. 
		
		Hence our model predicts that $ V_{\text{model}}(x,t) =   \frac{5.02365 \times 10^{-6}}{\log^{3.13788}(x)} \log^{4.3}(t) $. If we try plugging in $ x = \frac{t}{2 \pi} $, we get $ V_{\text{model}}(t) \sim 5.02365 \times 10^{-6} {\log^{1.16}(t)} $. This certainly means the errors for $ \zeta_{\frac{t}{2 \pi}} $ increase logarithmically for higher ordinates. However this model isn't quite right as it misses our plot by some margin:
		\newpage
		
		\begin{figure}[!h]
			
			\centering
			\includegraphics[scale=1]{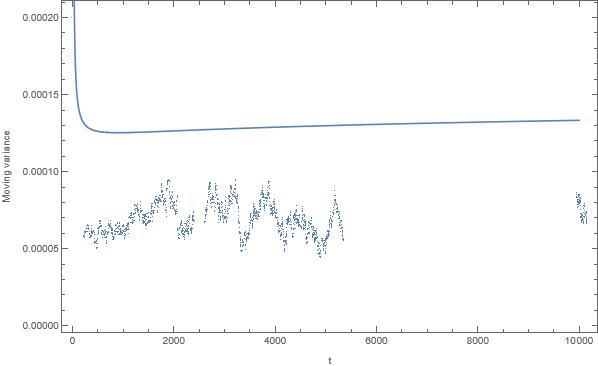}
			\caption{Model for $ B(t) = \log^{4.3}(t) $ and V(t)}
			\label{fig: logpow}
			
		\end{figure}
		We should try some other models for $ B(t) $ and compare.
		
	\end{enumerate}
	\subsection{Model $  B(t) = W(t)^k $}
	
	$ B(t) $ is modelled as $ W(t)^{3.7} $. An ideal $ k $ turns out to be something around 3.7. Finding a suitable $ A(x) = \frac{2.528 \times 10^{-5}}{W^{3.6}(x)} $, our model becomes:
	
	\newpage
	\begin{figure}[!h]
		
		\centering
		\includegraphics[scale=1]{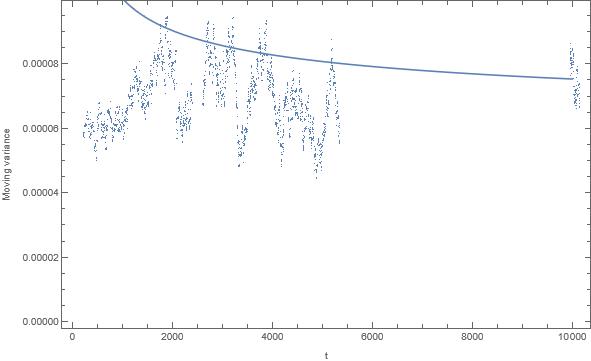}
		\caption{Model for $ B(t) = W(t)^{3.7} $ and V(t)}
		\label{fig: logpowAX}
		
	\end{figure}
	
	We see that the model is somewhat close to the actual values of the variance recorded. 
	
	\begin{comment}
		\subsection{Model $  B(t) = W(t)^{1.5} \log^{2.5}(t) $}
		
		$ B(t) $ is modelled as $W(t)^2 \log^2(t)  $. An ideal $ k $ turns out to be something around 3.7. Finding a suitable $ A(x) = \frac{1.37868 \times 10^{-5}}{\log ^3(x)} $, our model becomes:\\ 
		\begin{figure}[!h]
			
			\centering
			\includegraphics[scale=1]{modellogprodlog}
			\caption{Model for $ B(t) = W(t)^{3.7} $ and V(t)}
			\label{fig: logpowAX}
			
		\end{figure}\\
	\end{comment}
	
	\section{Where to go from here}
	\begin{enumerate}
		\item We assumed that the moving variances for fixed $ X $ at ordinate $ t $ will behave like  $ V(X,t) = A(X)B(t) $  However it may not be such a simple separable formula. 
		\item Look at even higher ordinates to see if the variances for $ X = t/{2 \pi} $ increase logarithmically or remain constant.
		\item Try to find a theoretical way to model the variances for fixed X, using Random Matrix Theory.
		\item Does a similar ``logarithmic" growth also appear in the case of general L Functions associate to elliptic curves. Assuming the errors remain constant, can it be used to calculate order of the zero at the central point. 
	\end{enumerate}

	\bibliographystyle{siam}
	\bibliography{References}
	
\end{document}